\documentclass[12pt]{article}
\usepackage{epsfig,latexsym,amsfonts,amsmath,amssymb}

\usepackage{colortbl}%
\usepackage[margin=3cm]{geometry}

\usepackage[notref,notcite]{showkeys}

\def\R{\mathbb{R}}

\def\f{\varphi}

\def\div{{\rm div}}

\def\irn{\int\limits_{\R^n}}

\def\DsD{\left(-\Delta\right)^{\!s}\!}
\def\DshalfD{\left(-\Delta\right)^{\!\frac{s}{2}}}

\def\DsN{\left(-\Delta_\Omega\right)^{\!s}}
\def\DshalfN{\left(-\Delta_\Omega\right)^{\!\frac{s}{2}}}

\def\DsNeps{\left(-\Delta_{\Omega_h}\right)^{\!s}\!}

\def\DsNepsh{(-\Delta_{\Omega_h})^{s}}

\def\DsNtilde{\left(-\Delta_{\widetilde\Omega}\right)^{\!s}\!}

\def\eps{\varepsilon}

\def\proof{\noindent{\textbf{Proof. }}}
\def\QED{\hfill {$\square$}\goodbreak \medskip}

\newtheorem{Theorem}{Theorem}[section]
\newtheorem{Lemma}[Theorem]{Lemma}

\newtheorem{Corollary}[Theorem]{Corollary}
\newtheorem{Remark}[Theorem]{Remark}
\newtheorem{Definition}[Theorem]{Definition}

\begin{document}

\title 
{Variational inequalities \\for the spectral fractional Laplacian}

\author{Roberta Musina\footnote{Dipartimento di Matematica ed Informatica, Universit\`a di Udine,
via delle Scienze, 206 -- 33100 Udine, Italy. Email: {roberta.musina@uniud.it}. 
{Partially supported by Miur-PRIN 201274FYK7\_004.}}~ and
Alexander I. Nazarov\footnote{
St.Petersburg Department of Steklov Institute, Fontanka 27, St.Petersburg, 191023, Russia, 
and St.Petersburg State University, 
Universitetskii pr. 28, St.Petersburg, 198504, Russia. E-mail: al.il.nazarov@gmail.com.
Supported by RFBR grant 14-01-00534.}
}
\date{}

\maketitle

\begin{abstract}
\footnotesize
In this paper we study the obstacle problems for the Navier (spectral) fractional Laplacian
$\DsN$ of order $s\in(0,1)$, in a bounded domain $\Omega\subset\R^n$.
\medskip

%
\end{abstract}

\normalsize

\bigskip\bigskip\bigskip

\section{Introduction}
Let $\Omega$ be a bounded and Lipschitz  domain in $\R^n$, $n\ge 1$. 
Given $s\in(0,1)$, a measurable function $\psi$ on $\Omega$ and $f\in \widetilde H^s(\Omega)'$, 
we consider the 
variational inequality
\begin{equation*}\tag{$\mathcal P_\Omega(\psi,f)$}
\label{eq:vi}
u\in K^s_\psi~,\qquad
\langle \DsN  u-f, v-u\rangle\ge 0\quad\forall v\in K^s_\psi~\!.
\end{equation*}
Here $\DsN$ is the {\em spectral} (or {\em Navier}) fractional
Laplacian, that is the $s$-th power of the standard Laplacian in the sense of spectral theory,
and
$$
K^s_\psi=\left\{v\in{\widetilde H^s(\Omega)}~|~v\ge \psi ~\textit{a.e. on ${{\Omega}}$}~\right\}~\!.
$$
We will always assume that the closed and convex set $K^s_\psi$
is not empty, also when not explicitly stated.

Problem \ref{eq:vi} admits a unique solution $u$, that can be characterized as
the unique minimizer for
\begin{equation}
\label{eq:minimization}
\inf_{v\in K^s_\psi}~\frac{1}{2}\langle \DsN  v,v\rangle-\langle f, v\rangle~\!.
\end{equation} 
The variational inequality \ref{eq:vi} is naturally related to 
the free boundary problem
\begin{equation}
{\label{eq:problem}}
\begin{cases}
u\ge \psi~,\quad \DsN u\ge f&\text{in ${\Omega}$}\\
\DsN u= f&\text{in $\{u>\psi\}$}\\
u=  0&\text{in  $\R^n\setminus\overline\Omega$}
\end{cases}
\end{equation}
as well. 
In fact, it is easy to show that any solution $u\in \widetilde H^s(\Omega)$ to (\ref{eq:problem})
satisfies \ref{eq:vi}, see Remark \ref{R:uniqueness}. The converse needs more care.
One of the main motivations of the present paper was indeed to find out
mild regularity assumptions on the data, to have that  the solution $u$
to \ref{eq:vi} solves the free boundary problem (\ref{eq:problem}).

Problem (\ref{eq:problem}) has been largely investigated
in case $\DsN$ is replaced by the {\em ''Dirichlet''} Laplacian $\DsD$, that is defined via the
Fourier transform by
\begin{equation}
\label{eq:fourier}
\mathcal F[\DsD u](\xi)=|\xi|^{2s}\mathcal F[u](\xi)= \frac{|\xi|^{2s}}{(2\pi)^{n/2}}\irn e^{-i~\!\!\xi\cdot x}u(x)~\!dx~\!,
\end{equation}
see Section \ref{S:ND} for details.
On this subject, we cite the pioneering paper  \cite{Sil} by Louis E. Silvestre, 
  \cite{BCRo, CaFi, CSS, GP, MNS, PP, SV.obs} and references there-in, with no attempt to 
provide a complete  list. 

As far as we know, the variational inequality
\ref{eq:vi} and the free boundary problem (\ref{eq:problem}) has never been discussed before. 
Actually, the {\em ''Navier'' case} is  more challenging, because of the dependence
of the differential operator $\DsD$ on the domain. This extra difficulty led 
us to investigate in Section \ref{S:ST} the dependence $\Omega\mapsto \DsN$.
The results there, and in particular Lemmata \ref{L:eige_ueps}, \ref{L:eige_ueps2}
might have an independent interest and could open new research directions,
see Remark \ref{R:Gamma}.

In Section \ref{S:truncations} we focus our attention on the action of $\DsN$ on truncations
$v\mapsto v^\pm$. The results here are essentially used in the remaining part
of the paper, as they provide the needed tools to construct test functions
for \ref{eq:vi}. Once we have developed the above mentioned tools, we
indicate how to modify the arguments in \cite{MNS} to find out  useful 
equivalent formulations and continuous dependence results for \ref{eq:vi},
see Section \ref{S:equivalent}.

Section \ref{S:regularity1} is entirely dedicated to regularity results. 
Most of the proofs here follows the outlines of the proofs in \cite{MNS}.
However here again more attention is needed because of the dependence
$\Omega\mapsto\DsN$; the preliminary results in Section \ref{S:ST}
will be crucially used in the proof of our main regularity result, that is
stated in Theorem \ref{T:measure}.

In the last section we take $f=0$ and compare the solution to \ref{eq:vi} with the solution
of the corresponding variational inequality with Dirichlet fractional operator $\DsD$.
The main result is stated in Theorem \ref{T:comparing1}.

%

\bigskip

\noindent
{\footnotesize
{\bf Notation.} For a bounded and Lipschitz domain  $\Omega\subset\R^n$ 
we denote by $-\Delta_\Omega$ the conventional Dirichlet Laplacian in $\Omega$, that is
the self-adjoint operator in $L^2(\Omega)$ defined by its quadratic 
form 
$$\langle-\Delta_\Omega u,u\rangle=\|\nabla u\|^2_2, \qquad u\in H^1_0(\Omega).
$$ 

We denote by $\lambda_j$, $j\ge 1$, the eigenvalues of $-\Delta_\Omega$ arranged in a non-decreasing unbounded sequence, according to their
multiplicities. Corresponding eigenfunctions
$$
\f_j\in H^1_0(\Omega)~,\qquad -\Delta\f_j=\lambda_j\f_j~,\quad  j\ge 1,
$$
form an orthogonal bases in $L^2(\Omega)$ and in $H^1_0(\Omega)$, and we assume them orthonormal in $L^2(\Omega)$. 

For $u\in H^1_0(\Omega)$ we have
\begin{equation*}
\label{eq:H_basis}
u=\sum_{j=1}^\infty \Big(\int\limits_\Omega u\f_j\Big)\f_j~,\quad -\Delta_\Omega u=\sum_{j=1}^\infty \lambda_j~\!\Big(\int\limits_\Omega u\f_j\Big)\f_j,
\end{equation*}
where the first series converges in $H^1_0(\Omega)$, while the second one has to be intended on the sense of distributions. Thus
$$
\|u\|_2^2=\sum_{j=1}^\infty \Big(\int\limits_\Omega u\f_j\Big)^2~,\quad \langle-\Delta_\Omega u,u\rangle= 
\sum_{j=1}^\infty \lambda_j\Big(\int\limits_\Omega u\f_j\Big)^2~\!.
$$

Next, take $s\in(0,1)$. The ``Navier'' (or spectral) fractional Laplacian of order $s$ on $\Omega$ is defined by the series
(in the sense of distributions) 
$$
\DsN\!u=\sum_{j=1}^\infty \lambda_j^s~\Big(\int\limits_\Omega u\f_j\Big)\f_j~\!.
$$
 It is known that the domain of the corresponding quadratic form $\langle\DsN u,u\rangle$ is the space
$$
{\widetilde H^s(\Omega)}:=\{u\in {H^s(\R^n)}~|~u\equiv 0~~\text{on}~~\R^n\setminus\overline\Omega\},
$$
see for instance \cite[Lemma 1]{FL}. The standard reference for the Sobolev space $H^s(\R^n)$ is the monograph
\cite{Tr} by Triebel.
  We endow $\widetilde H^s(\Omega)$ with
the Hilbertian norm
$$
\|u\|_{\widetilde H^s(\Omega)}^2
= \langle\DsN u,u\rangle=\|\DshalfN u\|_2^2~\!.
$$
Notice that $\f_j$ is the eigenfunction of $\DsN$ corresponding to the eigenvalue $\lambda_j^s$. That is, $\DsN$ is the $s$-th power of $-\Delta_\Omega$ 
in the sense of spectral theory.

We recall here some basic facts from \cite{ST}. For $u\in\widetilde H^s(\Omega)$ we put
\begin{gather*}
{\mathcal E}(w)=\!\!\int\limits_0^\infty\!\int\limits_{\R^n} y^{1-2s}|\nabla w(x,y)|^2\,dxdy~,\\
{\cal W}^\Omega_u=\Big\{w(x,y)~\!|~\!
{\mathcal E}(w)<\infty~,\ \ w\big|_{y=0}=u~,~w(\cdot, y)\equiv 0~~\text{on $\R^n\setminus\overline\Omega$}\Big\}.
\end{gather*}
The minimization problem
\begin{equation*}\tag{$\mathcal M^\Omega_u$}
\label{eq:CSmini}
\inf_{w\in{\cal W}^\Omega_u} {\mathcal E}(w)
\end{equation*}
has a unique solution $w^\Omega_u:\R^n\times\R_+~\longrightarrow~\R$, that solves the Dirichlet problem
\begin{equation*}\tag{$\mathcal{L}^\Omega_u$}
\label{eq:ST}
\begin{cases}
-\div (y^{1-2s}\nabla w)=0\qquad\text{in $\Omega\times\mathbb R_+$};\\
w(\cdot,y)\equiv 0\quad\text{on $(\R^n\setminus\overline\Omega)\times\R_+$},\quad
w\big|_{y=0}=u.
\end{cases}
\end{equation*}
The  results  in   \cite[Theorem 1.1]{ST} (see also Section 2 there-in), and integration by parts imply that
\begin{equation}
\label{eq:ST_energy}
\langle\DsN u,u\rangle=  {c_s}~\! {\mathcal E}(w^\Omega_u)~\!,
\end{equation}
for an explicitly known constant $c_s>0$. In addition, we have that
\begin{equation}
\label{eq:neumann}
\DsN u(x)=-c_s~\!\lim\limits_{y\to0^+} y^{1-2s}\partial_yw_u^\Omega(x,y)=
-2sc_s\lim\limits_{y\to0^+} \frac{w_u^\Omega(x,y)-u(x)}{y^{2s}},
\end{equation}
where the limits have to be intended in the sense of distributions.   

From (\ref{eq:ST}) and (\ref{eq:neumann}) it follows that for any function $w$ on $\R^n\times\R_+$ 
with finite energy $\mathcal E(w)$, one has
\begin{equation}
\label{eq:w_variational}
c_s \int\limits_0^\infty\!\int\limits_{\R^n} y^{1-2s}\nabla w^\Omega_u\cdot\nabla w\,dxdy=
\langle\DsN u, w\big|_{y=0}\rangle~\!.
\end{equation}}

\section{Dependence of the Navier Laplacian on $\Omega$}
\label{S:ST}

The following statement was in fact proved in \cite{FL}, see also \cite{FL2}. We give it with full proof for reader's convenience.

\begin{Lemma}
\label{L:lemma2}
Let $\Omega,\widetilde\Omega $ be bounded and Lipschitz domains in $\R^n$, with $\overline\Omega\subset\widetilde\Omega$.
Let $u\in \widetilde H^s(\Omega)$, $u\not\equiv 0$. Then
\begin{equation}
\label{eq:monotonicity}
\langle\DsNtilde u, u\rangle<\langle\DsN u, u\rangle.
\end{equation}
Moreover, if $u\ge0$ then $\DsNtilde u< \DsN u$ in the distributional sense on $\Omega$.
\end{Lemma}

\proof
Let $w^\Omega_u$, $w^{\widetilde\Omega}_u$ be the solutions to the minimization problems 
$(\mathcal M^{\Omega}_u)$, $(\mathcal M^{\widetilde\Omega}_u)$, respectively. 
Since $\widetilde H^s(\Omega)\subset \widetilde H^s(\widetilde\Omega)$, then  ${\cal W}^\Omega_u\subset {\cal W}^{\widetilde\Omega}_u$
and therefore ${\mathcal E}(w^{\widetilde\Omega}_u)\le{\mathcal E}(w^{\Omega}_u)$.
Now notice that $w^{\widetilde\Omega}_u$ is nontrivial and analytic on $\widetilde\Omega\times\R_+$,
because it solves ($\mathcal{L}^{\widetilde\Omega}_u$).  
Hence it can not vanish in $(\widetilde\Omega\setminus\overline\Omega)\times\R_+$, 
that trivially gives $w^{\widetilde\Omega}\neq  w^\Omega_u$ and ${\mathcal E}(w^\Omega_u)<{\mathcal E}(w^{\widetilde\Omega}_u)$.
Inequality (\ref{eq:monotonicity}) readily follows from
(\ref{eq:ST_energy}).

Notice that if $u\ge0$ then $w^{\widetilde\Omega}_u$ is positive on $\widetilde\Omega\times\mathbb R_+$ by the maximum principle. 
The function
$W:=w^{\widetilde\Omega}_u-w^{\Omega}_u$ solves
$$
\begin{cases}
-\div (y^{1-2s}\nabla W)=0\qquad\text{in $\Omega\times\mathbb R_+$}\\
W(\cdot,y)> 0\quad\text{on $(\widetilde\Omega\setminus\overline\Omega)\times\R_+$}~,\quad
W\big|_{y=0}=0,
\end{cases}
$$
and the maximum principle gives $W> 0$ on $\Omega\times\mathbb R_+$. Applying
the Hopf-Oleinik boundary point lemma (see \cite{KH}, \cite{ABMMZ}) to the function $W(x,t^{\frac 1{2s}})$, we obtain
$$
0< 2s\lim_{y\to 0^+}\frac{W(x,y)-W(x,0)}{y^{2s}}=
\lim_{y\to 0^+}y^{1-2s}\big(\partial_yw^{\widetilde\Omega}_u-
\partial_yw^{\Omega}_u\big)~,\quad x\in\Omega.
$$
The conclusion readily follows from (\ref{eq:neumann}).
\QED


Now let $\Omega_h$, $h\ge 1$, be a sequence of uniformly 
Lipschitz domains such that 
\begin{equation}
\label{domains}
\overline\Omega\subset \Omega_h\subseteq\Big\{ x\in\R^n~|~d(x,\Omega)<\frac1h\Big\}. 
\end{equation}
It is convenient to regard at $\widetilde H^s(\Omega)$, $s\in(0,1]$, 
as a subspace of $\widetilde H^s(B_R)$, where $B_R$ is an open  ball containing $\overline\Omega_h$, so that we 
have continuous embeddings
$\widetilde H^s(\Omega)\hookrightarrow \widetilde H^s(\Omega_h)\hookrightarrow \widetilde H^s(B_R)$. In particular,
 $H^1_0(\Omega)\hookrightarrow H^1_0(\Omega_h)\hookrightarrow H^1_0(B_R)$.

It is easy to show that the domains $\Omega_h$ $\gamma$-converge to $\Omega$ as $h\to\infty$. That is,
for any $f\in L^2(B_R)$, we have $v_{h}\to v$ in $H^1_0(B_R)$,
where the functions $v_{h}\in H^1_0({\Omega_h})$ and $v\in H^1_0(\Omega)$ are defined via
$$
-\Delta_{\Omega_h}\!v_{h}=f\quad\text{in ${\Omega_h}$}~,\quad -\Delta_{\Omega}v=f\quad\text{in $\Omega$.}
$$
Let us recall some facts from \cite[Lemma XI.9.5]{DS}, see also \cite[Theorem 2.3.2]{He} and \cite[Example 2.1]{BDM}. 
The $\gamma$-convergence of the domains 
implies that the eigenvalues and eigenfunctions of $-\Delta_{\Omega_h}$
converge, respectively, to the eigenvalues and eigenfunctions of $-\Delta_{\Omega}$.
More precisely, it turns out that
\begin{equation}
\label{eq:i)}
\lambda_j^h\to \lambda_j~,\qquad\f_j^h\to \f_j\quad\text{in $H^1_0(B_R)$\qquad as $h\to\infty$,}
\end{equation}
provided that eigenfunctions corresponding to multiple eigenvalues are suitably chosen.\medskip

Now we start to study the behavior of the fractional Laplacian $\DsNeps $ as $h\to \infty$.
The next lemma, of independent interest, will be crucially used in the proof of our regularity results.

\begin{Lemma}\label{L:eige_ueps}
Let $u_{h}\in \widetilde H^s(\Omega_h)$ be a bounded sequence
in $\widetilde H^s(\Omega_h)$ such that $u_{h}\to u$  in $L^2(B_R)$. Then
$u\in \widetilde H^s(\Omega)$ and 
$$
{\langle\DsN u,u\rangle\le \liminf_{h\to \infty}\langle\DsNepsh u_{h},u_{h}\rangle}~\!.
$$
\end{Lemma}

\proof
Clearly, $u$ is the weak limit of the sequence $u_h$ in $\widetilde H^s(B_R)$
and $u_h\to u$ almost everywhere. Hence $u\in\widetilde H^s(\Omega)$.  
Next, for any  integer $m\ge 1$ we have that
$$
\langle\DsNepsh u_h,u_h\rangle\ge  \sum_{j=1}^m (\lambda^h_j)^s\Big(\int\limits_\Omega u_h\f^h_j\Big)^2=
\sum_{j=1}^m \lambda_j^s\Big(\int\limits_\Omega u\f_j\Big)^2+o_h(1)
$$
by (\ref{eq:i)}). Thus 
$$
\liminf_{h\to\infty}\langle\DsNepsh u_h,u_h\rangle\ge \sum_{j=1}^m \lambda_j^s \Big(\int\limits_\Omega u\f_j\Big)^2~\!.
$$
Taking the limit as $m\to \infty$ we infer
$$
\liminf_{h\to \infty}\langle\DsNeps u,u\rangle\ge \sum_{j=1}^\infty \lambda_j^s\Big(\int\limits_\Omega u\f_j\Big)^2=
\langle\DsN u,u\rangle,
$$
that ends the proof. 
\QED

\begin{Lemma}\label{L:eige_ueps2}
Let $u\in \widetilde H^s(\Omega)$. Then
\begin{itemize}
 \item[$i)$] $\displaystyle{\lim_{h\to\infty}\langle\DsNeps u,u\rangle=\langle\DsN u,u\rangle}$;
 \item[$ii)$] $\DsNeps u\to\DsN u$ weakly in $\widetilde H^s(\Omega)'$~\!.
\end{itemize}
\end{Lemma}

\proof
The first claim is an immediate consequence of Lemmata \ref{L:lemma2} and \ref{L:eige_ueps}. Now take 
any test function ${v}\in \widetilde H^s(\Omega)$ and use $i)$  to get
$$
\begin{alignedat}[t]{2}
4\langle\DsNeps u,{v}\rangle&=\langle\DsNeps (u+{v}),(u+{v})\rangle-\langle\DsNeps (u-{v}),(u-{v})\rangle&\\
                &=\langle\DsN (u+{v}),(u+{v})\rangle-\langle\DsN (u-{v}),(u-{v})\rangle+o_h(1)&\\
                &=4\langle\DsN u,{v}\rangle+o_h(1).&\text{\hfill$\square$}
\end{alignedat}$$

\begin{Remark}
\label{R:Gamma}
Let us introduce the functionals $L^2(B_R)\to \R\cup\{\infty\}$,
$$
Q^{\Omega}_s(u)=\begin{cases}
\langle\DsN u,u\rangle&\text{if $u\in \widetilde H^s(\Omega)$;}\\
\infty&\text{otherwise.}
\end{cases}
$$
Lemma \ref{L:eige_ueps} and $i)$ in Lemma \ref{L:eige_ueps2} say that $Q^{\Omega}_s$ is the $\Gamma$-limit of the sequence $Q^{\Omega_h}_s$. One can 
wonder if this fact holds for any sequence of perturbating domains $\Omega_h$ that 
$\gamma$-converges to $\Omega$. 
For $s=1$, answer is positive, see \cite[Theorem 13.12]{DM}. 
Differently from the fractional case $s\in(0,1)$, for $s=1$ the quadratic forms $Q^{\Omega_h}_1$ and $Q^{\Omega}_1$ coincide on the intersection of
their domains.
\end{Remark}

\section{Truncations}
\label{S:truncations}
Truncation operators play an important role in studying  obstacle problems. For measurable functions $v,w$ we put
$$
v\vee w=\max\{v,w\}~,\quad v\wedge w=\min\{v,w\}~,\quad v^+=v\vee 0~,\quad v_-=-(v\wedge 0),
$$
so that $v=v^+-v^-$ and $|v|=v^++v^-$. It is well known that 
$v\vee w\in {H^s(\R^n)}$ and $v\wedge w\in {H^s(\R^n)}$
if $v,w\in {H^s(\R^n)}$. In addition it holds that
$$
(v+m)^-~,\quad  (v-m)^+~,\quad  v\wedge m\in {\widetilde H^s(\Omega)}
$$
for any $v\in {\widetilde H^s(\Omega)}$, $m\ge 0$, see  \cite[Lemma 2.4]{MNS}.

\begin{Lemma}
\label{L:m_new}
Let $v\in {\widetilde H^s(\Omega)}$ and  $m\ge 0$. 

$i)$ If $(v+m)$ changes sign, then
$${\langle {\DsN} v,(v+m)^-\rangle+\|\DshalfN (v+m)^-\|_2^2< 0};
$$

$ii)$ If $(v-m)$ changes sign, then 
$${\langle {\DsN} v,(v-m)^+\rangle-\|\DshalfN (v-m)^+\|_2^2> 0};
$$

$iii)$ If $(v-m)$ changes sign, then 
$${\|\DshalfN (v\wedge m)\|_2^2<\|\DshalfN v\|_2^2- \|\DshalfN (v- m)^+\|_2^2}.
$$
\end{Lemma}

\proof
Let $w^\Omega_v$ be the solution to the minimization problem
$(\mathcal M^\Omega_v)$, and let $w^\Omega_{(v+m)^-}$ be the 
solution to $(\mathcal M^\Omega_{(v+m)^-})$. Since ~$(w^\Omega_v+m)^-\in {\cal W}^\Omega_{(v+m)^-}$, we have   
\begin{equation}
\label{eq:CSstrict}
{\mathcal E}((w^\Omega_v+m)^-)\ge{\mathcal E}(w^\Omega_{(v+m)^-})~\!.
\end{equation}
In addition, from (\ref{eq:w_variational}) we get  
\begin{equation}
\label{eq:WW}
c_s \int\limits_0^\infty\!\int\limits_{\R^n} y^{1-2s}\nabla w^\Omega_v\cdot\nabla (w^\Omega_v+m)^-\,dxdy=
\langle\DsN v, (v+m)^-\rangle.
\end{equation}
It is well known that $\nabla w^\Omega_v\cdot\nabla (w^\Omega_v+m)^-=-|\nabla (w^\Omega_v+m)^-|^2$ a.e. in $\Omega\times\R_+$.
Thus (\ref{eq:WW}), (\ref{eq:CSstrict}) and (\ref{eq:ST_energy}) give
\begin{eqnarray}
 \label{eq:lemma}
 \nonumber
\langle\DsN v, (v+m)^-\rangle &=& -c_s{\mathcal E}((w^\Omega_v+m)^-)\\
&\le& -c_s{\mathcal E}(w^\Omega_{(v+m)^-}) \le -\|\DshalfN (v+m)^-\|_2^2~\!,
\end{eqnarray}
and $i)$ with a large inequality follows. 

Now assume that equality holds in $i)$. We have to show that $v+m$ is nonnegative or nonpositive.
Since equality holds everywhere in (\ref{eq:lemma}), then ${\mathcal E}((w^\Omega_v+m)^-)={\mathcal E}(w^\Omega_{(v+m)^-})$. 
We infer that 
$(w^\Omega_v+m)^-=w^\Omega_{(v+m)^-}$, as the minimization problem $(\mathcal M^\Omega_{(v+m)^-})$ admits a unique solution.
Whence, $(w^\Omega_v+m)^-$ solves ($\mathcal{L}^\Omega_{(v+m)^-}$) with nonnegative boundary datum $(v+m)^-$. By the maximum principle
either $(w^\Omega_v+m)^-\equiv0$, i.e. $v+m\ge0$; or $(w^\Omega_v+m)^->0$ in $\Omega\times\R_+$, that is,
$w^\Omega_v+m<0$ in $\Omega\times\R_+$ and $v+m\le 0$ on $\Omega$. This gives $i)$.  
To check $ii)$ notice that $(v-m)^+=((-v)+m)^-$ and then use $i)$ with $(-v)$ instead of $v$.

Finally, we write $v\wedge m=v-(v-m)^+$ and use $ii)$ to get
\begin{eqnarray*}
\|\DshalfN (v\wedge m)\|_2^2
&=&\|\DshalfN v\|_2^2-2
\langle\DsN v,(v-m)^+\rangle
+ \|\DshalfN (v-m)^+\|_2^2\\
&<&
\|\DshalfN v\|_2^2-
\|\DshalfN (v-m)^+\|_2^2~\!.
\end{eqnarray*}
Thus $iii)$ holds true, and the lemma is completely proved.
\QED

\begin{Remark}
For the Dirichlet fractional Laplacian the inequalities $i)$--$iii)$ were proved in \cite[Lemma 2.4]{MNS}
with large signs. Arguing as above and using the Caffarelli-Silvestre extension \cite{CaSi} instead of \cite{ST}
we can get complete ``Dirichlet'' analog of Lemma 4.1.
\end{Remark}

\begin{Remark}
\label{R:MP}
Taking $m=0$ in Lemma \ref{L:m_new} we obtain the ''Navier'' counterpart of \cite[Lemma 2.1]{MNS}.
The statement $iii)$ in this case can be rewritten as follows:
for $v\in {\widetilde H^s(\Omega)}$ with $v^+, v^-\neq 0$ 
$$\langle \DsN  v^+,v^-\rangle=\langle \DsN v^-,v^+\rangle< 0.
$$
Next, we notice that $ii)$ in Lemma \ref{L:m_new} with $m=0$ gives the well known
weak maximum principle for $\DsN$. A strong maximum principle was proved
in \cite[Lemma 2.4]{CDDS}. Namely, if $u\in \widetilde H^s(\Omega)\setminus\{0\}$
and $\DsN u\ge 0$ in $\Omega$ then $u$ is bounded away from zero on every compact set $K\subset\Omega$.
\end{Remark}

\begin{Remark}
\label{R:continuity}
Assume  $v_h\to v$ in $\widetilde H^s(\Omega)$. Then $\|\DshalfN (v_h^\pm-v)\|_2\to 0$. 
For the proof, recall that $\|\DshalfN\cdot\|_2$ equivalent to the norm in $\widetilde H^s(\Omega)$ that 
is induced by $H^s(\R^n)$,
see \cite[Corollary 1]{FL}; then use the continuity of truncation operators
$v\mapsto v^\pm$ in $H^s(\R^n)$, see for instance  \cite[Theorem 5.5.2/3]{RS}.
\end{Remark}

\section{Equivalent formulations\\ and continuous dependence results}
\label{S:equivalent}
We start by recalling the notion of (distributional) supersolution.
 
\begin{Definition}
A function ${\cal U}\in\widetilde H^s(\Omega)$ is a supersolution for $\DsN v=f$ if
$$
\langle \DsN{\cal U}-f,\f\rangle\ge  0\qquad\text{for any $\f\in \widetilde H^s(\Omega)$, $\f\ge 0$.}
$$
\end{Definition}
Thanks to the results in the previous section, the arguments in \cite[Section 3]{MNS} 
can be easily adapted to cover the problem \ref{eq:vi}. 
We start by pointing out some equivalent formulations for \ref{eq:vi}. For the proof, argue
as for \cite[Theorem 3.2]{MNS}.

\begin{Theorem} 
\label{T:sup}
Let $u\in K^s_\psi$. The following sentences are equivalent.
\begin{itemize}
\item[$a)$] $u$ is the solution to the problem \ref{eq:vi};
\item[$b)$] $u$ is the smallest supersolution for $\DsN v=f$ in the convex set $K^s_\psi$. That is, 
${\cal U}\ge u$ almost everywhere in ${\Omega}$,  for any supersolution ${\cal U}\in K^s_\psi$;
\item[$c)$]  $u$ is a supersolution for $\DsN v=f$ and 
$$\displaystyle{\langle\DsN u-f,(v-u)^-\rangle= 0}\qquad\text{for any $v\in K^s_\psi$.}
$$
\item[$d)$] $\displaystyle{\langle\DsN v-f,v-u\rangle\ge 0}$ for any $v\in K^s_\psi$.
\end{itemize}
\end{Theorem}

The next corollary is an immediate consequence of $a)\Rightarrow b)$ in Theorem \ref{T:sup}.

\begin{Corollary}\label{compare_f}
Let $f_1, f_2\in \widetilde H^s(\Omega)'$ and let
 $u_i$ be  the solution to $\mathcal P_\Omega(\psi,f_i)$, $i=1,2$.
If  $f_1\ge f_2$ in the sense of distributions, then $u_1\ge u_2$
a.e. in $\Omega$.
\end{Corollary}

\begin{Remark}
\label{R:uniqueness}
Let  $u\in K^s_\psi$ be such that $\DsN u\ge f$ in $\Omega$. Then
$\DsN u-f$ can be identified with a nonnegative Radon measure on $\Omega$.
Assume that the support of this measure is contained in the coincidence set $\{u=\psi\}$, 
so that $u$ solves the free boundary problem
(\ref{eq:problem}). Let $v\in K^s_\psi$. Since $(v-u)^-$ vanishes on $\{u=\psi\}$, we have
$\langle \DsN u-f,(v-u)^-\rangle=0$. Hence $u$ solves \ref{eq:vi} by Theorem \ref{T:sup}.
\end{Remark}

Now we can state our continuous dependence results. The proof of the
next theorem is totally similar to the proof of Theorem 4.1 in \cite{MNS}, and we omit it.

\begin{Theorem}
\label{T:bounded1}
Let $\psi_1,\psi_2$ be  given  obstacles, $f\in \widetilde H^s(\Omega)'$ and let
 $u_i$ be solutions to $\mathcal P_\Omega(\psi_i,f)$, $i=1,2$.
%
If $\psi_1-\psi_2\in L^\infty({\Omega})$, then the difference $u_1-u_2$ is bounded, and
$$
i)~~\|(u_1-u_2)^+\|_\infty\le\|(\psi_1-\psi_2)^+\|_\infty~,\quad
ii)~~\|(u_1-u_2)^-\|_\infty\le\|(\psi_1-\psi_2)^-\|_\infty.
$$
In particular, $\|u_1-u_1\|_\infty\le\|\psi_1-\psi_2\|_\infty$.
\end{Theorem}


\begin{Corollary}
\label{C:infty}
Let $\psi\in  L^\infty({\Omega})$ and $f\in L^p({\Omega})$, with $p\in(1,\infty)$, $p>n/2s$.
Let $u\in{\widetilde H^s(\Omega)}$ 
be the solution to \ref{eq:vi}.
Then $u\in L^\infty({\Omega})$ and
\begin{equation}
\label{eq:Linf}
\psi\vee\omega_f \le u\le \|\psi^+\|_\infty+c\|f^+\|_p\quad\text{a.e. in $\Omega$,}
\end{equation}
where $\omega_f$ solves the problem
\begin{equation}
\label{eq:e}
\DsN\omega_f=f\quad\text{in $\Omega$}~,\qquad \omega_f\in \widetilde H^s(\Omega),
\end{equation}
and $c$ depends only on $n,s,p$ and  ${\Omega}$.
In particular, if $f=0$ then 
$$
\psi^+\le u\le \|\psi^+\|_\infty~\!.
$$
\end{Corollary}

\proof 
Notice that $f\in \widetilde H^s(\Omega)'$ by Sobolev embedding theorem.
Since $u$ is supersolution of (\ref{eq:e}), the first inequality in (\ref{eq:Linf}) follows
by the maximum principle in Remark \ref{R:MP}. To prove the second inequality 
in (\ref{eq:Linf}) we introduce the functions 
$\omega^N_{f^+}, \omega^D_{f^+}\in \widetilde H^s(\Omega)$ via
$$
\DsN \omega^N_{f^+}=\DsD \omega^D_{f^+}=f^+\qquad\text{in $\Omega$}.
$$
It has been proved in \cite{MNS}, proof of Corollary 4.2, that $\omega^D_{f^+}\le c\|f^+\|_p$, where
the constant $c>0$  does not depend on $f$. 

Next, $\omega^D_{f^+}\ge 0$ by the maximum principle. Therefore, $\DsN \omega^D_{f^+}\ge \DsD \omega^D_{f^+}$ in 
$\Omega$ by \cite[Theorem 1]{FL}. Thus $\DsN(\omega^D_{f^+}-\omega^N_{f^+})\ge 0$
in $\Omega$, that implies $\omega^D_{f^+}\ge \omega^N_{f^+}$ by the maximum principle in Remark \ref{R:MP}. 
In particular we have $\omega^N_{f^+}\le c\|f^+\|_p$ a.e. in $\Omega$.

Now let $u_1$ be the unique solution of $\mathcal P_\Omega(\psi,f^+)$. 
We can consider $\omega^N_{f^+}$ as the solution of the problem $\mathcal P_\Omega(\omega_{f^+},f^+)$,
so that 
Theorem \ref{T:bounded1} gives
$$
u\le (u_1-\omega^N_{f^+})^++\omega^N_{f^+}\le \|(\psi-\omega^N_{f^+})^+\|_\infty+\omega^N_{f^+},
$$
and the second inequality in (\ref{eq:Linf}) readily follows.
\QED

To prove the next continuous dependence results, argue as  in \cite{MNS},
proofs of Theorems 4.3 and 4.4, respectively.

\begin{Theorem}
\label{T:Linfty}
Let $\psi_h\in L^\infty(\Omega)$ be  a sequence of obstacles and let $f\in \widetilde H^s(\Omega)'$ be given.
Assume that there exists $v_0\in \widetilde H^s(\Omega)$, such that $v_0\ge \psi_h$ for any $h$.

Denote by  $u_h$ the solution to the obstacle problem $\mathcal P_\Omega(\psi_h,f)$.
If $\psi_h\to \psi$ in $L^\infty(\Omega)$, then 
$u_h\to u$ in $\widetilde H^s(\Omega)$, where $u$ is the solution to the limiting problem
\ref{eq:vi}.
\end{Theorem}

 \begin{Theorem}
 \label{T:Hs2}
Let $\psi_h\in {H^s(\R^n)}$ be  a sequence of obstacles such that $\psi_h^+\in \widetilde H^s(\Omega)$,
and let $f_h$ be a sequence in 
$\widetilde H^s(\Omega)'$.
Assume that
$$
\psi_h\to \psi\quad\text{in $H^s(\R^n)$, and}\quad f_h\to f\quad\text{in $H^s(\R^n)'$}.
$$
Denote by  $u_h$ the solution to the obstacle problem
$\mathcal P_\Omega(\psi_h,f_h)$. 
Then $u_h\to u$ in ${\widetilde H^s(\Omega)}$, where $u$ is the solution of the limiting obstacle problem \ref{eq:vi}.
\end{Theorem}

\section{Regularity results}
\label{S:regularity1}

Let $u$ be the solution to \ref{eq:vi}. 
In this Section we provide estimates of the Radon measure $\DsN u-f\ge 0$ in $\Omega$ and
a regularity result for $u$. 
The results in Section \ref{S:ST} will be largely used.

Recall that $\omega_f$ is the solution to the boundary value problem (\ref{eq:e}).

\begin{Theorem}
\label{T:measure}
Assume that $f\in\widetilde H^s(\Omega)'$ and that $f, \psi$ satisfy the following conditions.
\begin{itemize}
\item[$A1)$] $(\psi-\omega_f)^+\in \widetilde H^s(\Omega)$;
\item[$A2)$] $\DsN(\psi-\omega_f)^+-f$ is a  locally finite signed measure on $\Omega$.
\end{itemize}
Let $u\in{\widetilde H^s(\Omega)}$ be the solution to \ref{eq:vi}.
Then 
$$
0\le \DsN u-f\le (\DsN(\psi-\omega_f)^+-f)^+\quad\text{in the distributional sense on $\Omega$.}
$$
\end{Theorem}

\proof
The first step is quite similar to the Dirichlet case (\cite[Theorem 1.1]{MNS})
and is based on the penalty method by Lewy-Stampacchia \cite{LS1}.

Notice that we can assume  $f=0$, 
$\psi\in \widetilde H^s(\Omega)$ and $\psi\ge 0$ in $\Omega$. In fact,
since
$u-\omega_f\in \widetilde H^s(\Omega)$ and $\DsN(u-\omega_f)\ge 0$,
then $u-\omega_f\ge 0$ in $\Omega$ by the maximum principle in Remark \ref{R:MP}. 
Thus $u-\omega_f\ge \psi\vee\omega_f-\omega_f=(\psi-\omega_f)^+$ and
$u-\omega_f$ solves the obstacle problem
$\mathcal P_\Omega((\psi-\omega_f)^+,0)$. 

In conclusion, we only have to show that 
\begin{equation}
\label{eq:ineq}
0\le \DsN u\le (\DsN\psi)^+\quad\textit{in the distributional sense on $\Omega$,}
\end{equation}
where $u$ solves $\mathcal P_\Omega(\psi,0)$ and $\psi$ is a nonnegative obstacle in $\widetilde H^s(\Omega)$, 
such that $\DsN\psi$
is a measure on $\Omega$. The first inequality in (\ref{eq:ineq}) holds by Theorem \ref{T:sup}. 

Now we prove the following claim:
\begin{equation}
\label{eq:step1}
\text{\em Assume ${\DsN\psi}\in L^p(\Omega)$ for any $p>1$. Then (\ref{eq:ineq}) holds.}
\end{equation}
We take $p\ge \frac{2n}{n+2s}$, that is needed only if $n>2s$. Then
$\widetilde H^s(\Omega)\hookrightarrow L^{p'}(\Omega)$ and
$L^{p}(\Omega)\subset \widetilde H^s(\Omega)'$ by Sobolev embeddings. In
particular $({\DsN\psi})^+\in \widetilde H^s(\Omega)'$. 
Take a function $\theta_\eps\in C^\infty(\R)$ such that $0\le \theta_\eps\le 1$, and
$\theta_\eps(t)=1$ for $t\le 0$, $\theta_\eps(t)=0$ for $t\ge \eps$. Let $u_\eps\in \widetilde H^s(\Omega)$
be the unique solution to 
\begin{equation}
\label{eq:penalty}
\DsN u_\eps=\theta_\eps(u_\eps-\psi)~\!({\DsN\psi})^+\quad\text{in $\Omega$.}
\end{equation}
Notice that $(1-\theta_\eps(u_\eps-\psi))(\psi-u_\eps)^+=0$ a.e. in $\Omega$. Therefore,
using $ii)$ in Lemma \ref{L:m_new} with $v=\psi-u_\eps$ and $m=0$ one gets $(\psi-u_\eps)^+\equiv 0$. 
In particular we infer that $u_\eps\in K^s_\psi$. On the other hand, $\DsN u_\eps\ge 0$; thus $u_\eps\ge u$
by $b)$ in Theorem \ref{T:sup}. 

Next, notice that $\DsN(u_\eps-u)\le \DsN u_\eps$ and that $\theta_\eps(u_\eps-\psi)(u_\eps-u-\eps)^+=0$ a.e. in $\Omega$.
Then again $ii)$ in Lemma \ref{L:m_new} plainly implies $(u_\eps-u-\eps)^+\equiv 0$. In conclusion, we have $u\le u_\eps\le u+\eps$, 
hence $\|u_\eps-u\|_\infty\to 0$ as $h\to\infty$. Therefore, for any nonnegative test
function $\f\in C^\infty_0(\Omega)$ we have that
\begin{eqnarray*}
\langle \DsN u,\f\rangle&=&\int\limits_\Omega u\DsN\f~\!dx=
\int\limits_\Omega u_\eps\DsN\f ~\!dx+o_\eps(1)\\
&=&\langle \DsN u_\eps,\f\rangle+o_\eps(1)
\le \langle (\DsN \psi)^+,\f\rangle+o_\eps(1).
\end{eqnarray*}
Thus, $\DsN u\le (\DsN \psi)^+$ in the distributional sense in $\Omega$, and (\ref{eq:step1}) is proved.

\medskip
The second step uses an approximation argument that requires to enlarge the domain $\Omega$.
It needs more care than in the Dirichlet case, because of the dependence of the Navier quadratic form on the domain.
Let $\Omega_h\supset\overline\Omega$, $h\ge 1$, be a sequence of uniformly Lipschitz domains satisfying (\ref{domains}).
The convex set
$$
K_h(\psi)=\{ v\in \widetilde H^s(\Omega_h)~|~v\ge \psi~~\text{a.e. on $\R^n$}~\}
$$
contains $K^s_\psi$, hence it is not empty. Let $u_h\in \widetilde H^s(\Omega_h)$ be the solution
to
\begin{equation*}\tag{$\mathcal P_{\Omega_h}$}
\label{eq:eps}
u_h\in K_h(\psi)~,\qquad
\langle \DsNeps  u_h, v-u_h\rangle\ge 0\qquad\forall v\in K_h(\psi)~\!.
\end{equation*}
We claim that 
\begin{equation}
\label{eq:claim1}
\DsNeps  u_h\le (\DsN \psi)^+\quad \textit{in the distributional sense on $\Omega$.}
\end{equation}
Fix $h$, and approximate $\psi$ with a sequence of smooth obstacles
$\psi^{k}=\psi*\rho_{k}$, where supp$(\rho_k)\subset B_{\frac{1}{k}}$. For $k$ large enough
we have $\psi^{k}\in C^\infty_0(\Omega_h)$. In addition
$\psi^{k}\to \psi$ in $\widetilde H^s(B_R)$ as $k\to\infty$, where $B_R$ is any ball containing $\Omega_h$. 
Now, let $u_{h}^k\in \widetilde H^s(\Omega_h)$ be the solution to the obstacle problem
\begin{equation*}\tag{$\mathcal P_{\Omega_h}^k$}
\label{eq:vih}
u_{h}^k\in K_{h}(\psi^{k})~,\qquad
\langle \DsNeps  u_{h}^k, v-u_{h}^k\rangle\ge 0 \qquad\forall v\in K_{h}(\psi^{k})~\!.
\end{equation*}
Then $u_{h}^k\to u_h$  in $\widetilde H^s(\Omega_h)$ as $k\to \infty$ by Theorem \ref{T:Hs2}, and (\ref{eq:step1}) gives
\begin{equation}
\label{eq:eps_h}
\DsNeps  u_{h}^k\le (\DsNeps \psi^{k})^+\quad\text{in the distributional sense on $\Omega$.}
\end{equation}
Next, $(\DsNeps \psi)^+*\rho_k$ is a nonnegative smooth function, and
$$(\DsNeps \psi)^+*\rho_k\ge (\DsNeps \psi)*\rho_k=\DsNeps\psi^{k}.$$
Thus
$(\DsNeps \psi)^+*\rho_h\ge (\DsNeps\psi^{k})^+$, and (\ref{eq:eps_h}) implies
$$
\DsNeps  u_h^k\le (\DsNeps \psi)^+*\rho_k\quad\text{in the distributional sense on $\Omega$.}
$$
Now, as $k\to \infty$ we have that
 $(\DsNeps \psi)^+*\rho_k\to (\DsNeps \psi)^+$ in the sense of measures,
and $\DsNeps u^k_h\to \DsNeps u_h$ in the sense of distributions. Thus
$$
\DsNeps  u_h\le (\DsNeps \psi)^+\quad\text{in the distributional sense on $\Omega$.}
$$
Since $\psi\in\widetilde H^s(\Omega)$ is nonnegative, Lemma \ref{L:lemma2}
gives {$(\DsNeps \psi)^+\le (\DsN \psi)^+$}, and  (\ref{eq:claim1}) follows.
\medskip

The last step is the passage to the limit in (\ref{eq:claim1}) as the domains $\Omega_h$ shrink to $\Omega$.
It makes the main difference with respect to the Dirichlet case.
We notice that
$u\in \widetilde H^s(\Omega_h)$ and in particular $u\in K_h(\psi)$. 
Therefore, using the variational characterization of $u_h$ as the solution to $(\mathcal P_{\Omega_h})$ 
and Lemma \ref{L:lemma2} we find
\begin{equation}
\label{eq:uff}
\langle  \DsNeps  u_h,u_h\rangle\le 
\langle  \DsNeps  u,u\rangle\le \langle  \DsN  u,u\rangle~\!.
\end{equation}
Lemma \ref{L:lemma2} gives also 
$\langle (-\Delta_{B_R}\!)^{s}u_h,u_h\rangle \le \langle  \DsNeps  u_h,u_h\rangle$.
Thus
(\ref{eq:uff}) implies that the sequence $u_h$ is bounded in $\widetilde H^s(B_R)$,
and therefore 
we can assume that $u_h\to \tilde u$ weakly in $\widetilde H^s(B_R)$. Using 
Lemma \ref{L:eige_ueps} and  (\ref{eq:uff}) we readily get $\tilde u\in \widetilde H^s(\Omega)$ and
\begin{eqnarray}
\langle\DsN \tilde u,\tilde u\rangle\le \liminf_{h\to\infty}\langle\DsNepsh u_h,u_h\rangle
\le  \limsup_{h\to\infty}\langle\DsNepsh u_h,u_h\rangle
\le\langle\DsN u,u\rangle,
\label{eq:party}
\end{eqnarray}
that is, $\langle  \DsN  \tilde u,\tilde u\rangle\le \langle  \DsN  u,u\rangle$.
On the other hand, $u_h\to \tilde u$ almost everywhere and 
$ u_h\ge \psi$ on ${\Omega}$. Thus $\tilde u\in K^s_\psi$.
Using the characterization of $u$ as the unique solution to the minimization problem (\ref{eq:minimization})
(with $f\equiv 0$), we first get $\tilde u=u$. Then we use (\ref{eq:party})
to infer that $u_h\to u$ in $\widetilde H^s(B_R)$. 

To conclude take any nonnegative function $\eta\in C^\infty_0(\Omega)$. 
Using $\overline\Omega\subset\Omega_h\subset\overline\Omega_h\subset B_R$ and Lemma \ref{L:lemma2}, we see that 
$(-\Delta_{B_R}~\!)^{s}\eta\le \DsNeps\eta\le \DsN\eta$ in $B_R$~\!.
Thus $\DsNeps\eta$ is a bounded sequence in $L^2(B_R)\hookrightarrow\widetilde H^s(B_R)'$. 
In addition, recall that
$u_h\to u$ in $\widetilde H^s(B_R)$ and $\DsNeps\eta\to \DsN\eta$ weakly in $\widetilde H^s(\Omega)'$,
see Lemma \ref{L:eige_ueps2}. Thus, we can use (\ref{eq:claim1}) to estimate
\begin{eqnarray*}
\langle(\DsN \psi)^+,\eta\rangle &\ge&
\langle\DsNeps u_h,\eta\rangle 
=\langle\DsNeps \eta,u_h\rangle
=\langle\DsNeps\eta,u\rangle+o_h(1)\\
&=&\langle\DsN \eta,u\rangle+o_h(1)= \langle\DsN u, \eta\rangle+o_h(1)~\!.
\end{eqnarray*}
The proof  is complete.
\QED


In the next results we adopt ''pointwise'' definitions of the contact set and of the non-contact set (compare with Definition \ref{Def} below for
a different notion), that is,
\begin{equation}
\label{eq:set}
\{u=\psi\}:=\{x\in\Omega~|~u(x)=\psi(x)\}~,\quad
\{u>\psi\}:=\{x\in\Omega~|~u(x)>\psi(x)\}~\!.
\end{equation}
Clearly, $\{u=\psi\}$ and $\{u>\psi\}$ are determined up to negligible sets. 

\begin{Theorem}
\label{T:regularity}
Let $\psi$ and $f$ satisfy assumptions of Theorem \ref{T:measure} and
\begin{itemize}
\item[$A3)$] $(\DsN(\psi-\omega_f)^+-f)^+\in L^p_{\rm loc}(\Omega)$ for some $p\in[1,\infty]$.
\end{itemize}
Let $u\in{\widetilde H^s(\Omega)}$ be the solution to \ref{eq:vi}.
Then the following facts hold.
\begin{itemize}
\item[$i)$] $\DsN u-f\in L^p_{\rm loc}({\Omega})$;
\item[$ii)$] $0\le \DsN u-f\le (\DsN(\psi-\omega_f)^+-f)^+$ a.e. on $\Omega$;
\item[$iii)$] $\DsN  u=f$ a.e. on $\{u>\psi\}$.
\end{itemize}
In particular, $u$ solves the free boundary problem (\ref{eq:problem}).
\end{Theorem}

\proof 
We follow the proof of \cite[Theorem 1.1]{MNS}.
Statements $i)$ and $ii)$ hold by Theorem \ref{T:measure}. Let us prove the last claim. As before, we assume 
$f\equiv 0$. Then $\DsN u\ge 0$. If $u\equiv0$ then $iii)$ is evident. Otherwise the strong maximum principle, see Remark \ref{R:MP}, gives $u>0$ in $\Omega$.
In particular, $u\ge \psi^+$ and $\{u>\psi\}=\{u>\psi^+\}$.

Use $c)$ in Theorem \ref{T:sup} with $v= \psi^+\in \widetilde H^s(\Omega)$, to get $\langle\DsN u,u-\psi^+\rangle =0$. 
Let $\f\in C^\infty_0(\Omega)$ be any nonnegative cut off function; for $m\ge 1$ put $g_m=(u-\psi^+)\wedge m$. 
Since $\DsN u\ge 0$ and $u-\psi^+\ge \f g_m$, we have
$$
0=\langle \DsN u,u-\psi^+\rangle\ge\langle \DsN u,\f g_m\rangle=\int\limits_{\Omega} {\DsN u}\cdot(\f g_m)dx.
$$
The last equality holds because $\DsN u\in L^1_{\rm loc}(\Omega)$ and  $\f g_m\in L^\infty(\Omega)$ has compact support in $\Omega$.
Thanks to the monotone convergence theorem we get
$$
0\ge \lim_{m\to\infty}\int\limits_{\Omega} {\DsN u}~\!\!\cdot \f g_m~\!dx=
\int\limits_{\Omega} ({\DsN u}~\!\!\cdot (u-\psi^+))\f~\!dx.
$$
Now $({\DsN u}~\!\!\cdot (u-\psi^+))\f\ge 0$ a.e. in $\Omega$, that gives
$({\DsN u}~\!\!\cdot (u-\psi^+))\f=0$ a.e. in $\Omega$. Since $\f$ was arbitrarily chosen, we 
conclude that $\DsN u~\!\!\cdot (u-\psi^+)=0$  a.e. in $\Omega$, and $iii)$ is proved.
\QED

\begin{Remark}
\label{R:regularity}
To obtain better regularity results for $u$, one can apply the regularity
theory for
$$
\DsN u=g\in L^p(\Omega)\quad\text{in $\Omega$}~,\qquad u\in \widetilde H^s(\Omega).
$$
In particular, if $p>\frac{n}{2s}$, then $u$ is H\"older continuous in $\overline\Omega$, see \cite[Corollary 3.5]{Gr}.
\end{Remark}

We conclude this section by giving a sufficient condition for the continuity of $u$.

\begin{Theorem}
\label{T:2}
Let $\psi\in C^0(\overline\Omega)$ be a given obstacle, such that $K^s_\psi$ is not empty and $\psi\le 0$ on $\partial\Omega$.
Let $f\in L^p(\Omega)$, with $p>n/2s$.
Then $u$ is continuous on $\R^n$.
\end{Theorem}

\proof
The argument is the same as in \cite[Theorem 1.2]{MNS}.
Fix a small $\eps>0$. We can assume that $\psi-\eps\in C^0_0(\R^n)$ and 
$\psi-\eps\le 0$ outside $\Omega$. Let $\psi^\eps_{h}$ be a sequence in $C^\infty_0(\R^n)$ such that
$\psi^\eps_h\to \psi-\eps$ uniformly on $\R^n$, as $h\to\infty$.

By Theorem \ref{T:regularity}, the solution $u_h\in \widetilde H^s(\Omega)$ to 
$\mathcal P_\Omega(\psi^\eps_h,f)$ satisfies $\DsN u^\eps_h\in L^p(\Omega)$ and therefore 
$u^\eps_h$ is H\"older continuous, see Remark \ref{R:regularity}.
Moreover, the estimates in Theorem \ref{T:bounded1} imply
that $u^\eps_h \to u^\eps$ uniformly on  $\Omega$, where $u^\eps$ solves
$\mathcal P_\Omega(\psi-\eps,f)$. In particular, $u^\eps\in C^0(\overline\Omega)$.
Finally, use again Theorem \ref{T:bounded1} to get that $u^\eps\to u$ uniformly, 
and conclude the proof.
\QED

\section{Comparing the Navier and the Dirichlet problems}
\label{S:ND}

In this section we compare the unique solutions $u_N, u_D$ to
the variational inequalities 
\begin{equation*}\tag{$\mathcal P_N(\psi)$}
\label{eq:N}
u_N\in K^s_\psi~,\qquad
\langle \DsN  u_N, v-u_N\rangle\ge 0\quad\forall v\in K^s_\psi
\end{equation*}
\begin{equation*}\tag{$\mathcal P_D(\psi)$}
\label{eq:D}
u_D\in K^s_\psi~,\qquad
~~\!\langle \DsD  u_D, ~\!v-u_D~\!\rangle~\ge 0\quad\forall v\in K^s_\psi~\!,
\end{equation*}
respectively. Here $\DsD$ is
the ''Dirichlet'' (or restricted) Laplacian, that has been already introduced in  (\ref{eq:fourier}). 

Problem \ref{eq:D}
has been investigated in \cite{MNS}. Recall that
\begin{equation}
\label{eq:L_positive}
\DsN u_N\ge 0~,\quad\DsD u_D\ge 0\qquad\text{on $\Omega$}
\end{equation}
in the sense of distributions, see Theorem \ref{T:sup} and \cite[Theorem 3.2]{MNS}, so that
$u_N, u_D$ are nonnegative by the maximum principle, and $u_N, u_D\ge \psi^+$. 
Thus, $u_N, u_D$ solve problems $\mathcal P_N(\psi^+)$, $\mathcal P_D(\psi^+)$, respectively. 
Hence we can assume without loss of generality that $\psi\ge0$ a.e. in $\Omega$.
Since $\psi\equiv0$ easily implies $u_N\equiv u_D\equiv0$, we assume also that $\psi\not\equiv0$.
By the strong maximum principle (see, respectively, \cite[Lemma 2.4]{CDDS} and 
\cite[Theorem 2.5]{IMS}) we have
\begin{equation}
\label{eq:u_positive}
u_N> 0\quad\text{and}\quad u_D> 0\qquad\text{on $\Omega$},
\end{equation}
in the sense that $u_N, u_D$ are bounded away from $0$ on every
compact set $K\subset\Omega$.

In this section we need to refine the notion of contact set introduced
in (\ref{eq:set}). We essentially use an idea 
due to Lewy and Stampacchia \cite{LS1, LS2}, see also \cite[Section 6]{KS}. 
We start with some preliminaries.

\begin{Definition}\label{Def}
Let $v$ be a nonnegative measurable function on the open set $\Omega$, and let $x\in \Omega$.
We say that $v(x)>0$ if
there exist $\rho,\eps>0$ such that $B_\rho(x)\subseteq\Omega$ and
$$
v-\eps\ge 0~~\text{almost everywhere in $B_\rho(x)$}~\!.
$$
\end{Definition}
For any measurable function $v$ on $\Omega$ we define
$$
P[v]=\{~\!x\in \Omega~|~v(x)>0~\!\}
$$
and we put $I[v]=\Omega\setminus P[v]$. The set $P[v]$ is clearly open; thus $I[v]$ is closed in $\Omega$. 
\begin{Lemma}
\label{L:meas1}
Let $v\ge 0$ be a  measurable function on $\Omega$ and let $K\subset P[v]$ be a compact set.
Then there exists $\eps_0>0$ such that $v-\eps_0\ge 0$ a.e. on $K$.

 In particular, $v>0$ a.e. in $P[v]$.
\end{Lemma}

\proof
If $K$ is not empty, for any $x\in K$ there exists a ball $B_x$ about $x$ such that
${B}_x\subset\Omega$ and $v\ge \eps_x>0$ on $B_x$. Since $K$ is compact,
we can find a finite number of points $x_i\in K$ such that $K$ is covered by the finite family $B_{x_i}$. Let $\eps_0=\min_i\eps_{x_i}$.
Then $v\ge \eps_0$ a.e. on $K$ and the first claim is proved. 

Now put $N=\{x\in P[v]~|~ v(x)=0\}$. By the first part of the proof, any compact set $K\subset N$ must have 
null measure. Thus $N$ is a negligile set.
\QED
By Lemma \ref{L:meas1} we have  the inclusions
$$
P[u_N-\psi]\subseteq\{u_N>\psi\}~,\quad P[u_D-\psi]\subseteq\{u_D>\psi\}.
$$
It might happen that 
$\{u_N>\psi\}$ and $\{u_D>\psi\}$ have positive measure but $P[u_N-\psi]$ and
$P[u_D-\psi]$ are empty, see Remark \ref{R:example} below.

\begin{Theorem}
\label{T:comparing1}
The following facts hold true.
\begin{itemize}
\item[$i)$] $\DsN u_N=0$ on $P[u_N-\psi]$ and $\DsD u_D=0$ on $P[u_D-\psi]$;
\item[$ii)$] $\DsN u_D> 0$ in the distributional sense on $\Omega$;
\item[$iii)$] $u_N\le u_D$; 
\item[$iv)$] $u_N< u_D$ on $P[u_N-\psi]$;
\item[$v)$] $u_D$ is the solution to the obstacle problem $\mathcal P_D(u_N)$;
\item[$vi)$] $\|\DshalfD u_D\|_2\le \|\DshalfD u_N\|_2<\|\DshalfN u_N\|_2\le \|\DshalfN u_D\|_2$,
and all signs are strict, unless $u_N\equiv u_D$.
\end{itemize}
\end{Theorem}

\proof
Fix any nonnegative $\f\in C^\infty_0(P[u_N-\psi])$. By Lemma \ref{L:meas1} we have that
$u_N\mp t\f\in K_\psi^s$ for sufficiently small $t\in\R$. Thus $\pm t\langle \DsN  u_N, \f\rangle\ge 0$, 
that proves $i)$ for $u_N$. The argument for $u_D$ is the same.

To prove $ii)$ use \cite[Theorem 1]{FL}, that gives $\DsN u_D>\DsD u_D\ge0$ by (\ref{eq:L_positive}).

By $b)$ in Theorem \ref{T:sup}, we know that $u_N$ is the smallest supersolution to $\DsN v=0$ 
in the set $K^s_\psi$. Thus $u_N\le u_D$ by $ii)$, and $iii)$ is proved. 

On the open set $P[u_N-\psi]$ both $u_N$ and $u_D$ are smooth by $i)$.
Assume by contradiction that there exists $x\in P[u_N-\psi]$ such that $(u_D-u_N)(x)=0$.
Then from $iii)$ we see that 
\begin{eqnarray*}
\DsD (u_D-u_N)(x)&=&C\cdot V.P. \int\limits_{\R^n}\frac {(u_D-u_N)(x)-(u_D-u_N)(y)}{|x-y|^{n+2s}}\,dy\\
&=&C\cdot V.P. \int\limits_{\R^n}\frac {-(u_D-u_N)(y)}{|x-y|^{n+2s}}\,dy \le 0.
\end{eqnarray*}
But this is impossible, as $\DsD(u_D-u_N)>\DsD u_D-\DsN u_N =0$ on $P[u_N-\psi]$ by \cite[Theorem 1]{FL} and $i)$.
Claim $iv)$ is  proved.

To check $v)$ we use \cite[Theorem 3.2]{MNS}, that characterizes 
$u_D$ as the smallest supersolution to $\DsD u=0$ in $K^s_\psi$. Since $K^s_{u_N}\subseteq K^s_\psi$, and $u_D\in K^s_{u_N}$ by $iii)$,
we see that $u_D$ is the smallest supersolution to $\DsD u=0$ in $K^s_{u_N}$ and the conclusion follows
again by \cite[Theorem 3.2]{MNS}.

Finally, Theorem 2 in \cite{FL} gives $\|\DshalfD u_N\|_2<\|\DshalfN u_N\|_2$.
The remaining inequalities in $vi)$ follow from the the fact that $u_N$ and $u_D$ are unique solutions to
\ref{eq:N} and \ref{eq:D}, respectively, and from variational formulations of these problems,
see (\ref{eq:minimization}) and \cite[(1.2)]{MNS} (with $f=0$).
\QED

\begin{Remark}
\label{R:example}
Take a smooth function $\eta\in \widetilde H^s(\Omega)$ satisfying
$\DsD \eta\ge 0$, $\eta>0$ in $\Omega$. Let ${\kappa}\subset\Omega$ be a compact set, having positive measure but empty interior. Consider
the obstacle $\psi=\eta\chi_{\Omega\setminus\kappa}$ and the solutions $u_N, u_D$ to \ref{eq:N}, \ref{eq:D}, respectively.
Clearly $\eta\in K^s_{\psi}$. Thus $\eta\ge u_D$ because
$u_D$ is the smallest supersolution for $\DsD v=0$ in $K^s_{\psi}$. But then
we have
$\eta\ge u_D\ge u_N\ge \psi=\eta\chi_{\Omega\setminus\kappa}$. In particular, $u_D= u_N={\psi}$ on $\Omega\setminus {\kappa}$.
Actually $\{u_D={\psi}\}=\{u_N={\psi}\}=\Omega\setminus {\kappa}$, because $u_N, u_D$ are positive in $\Omega$.
In particular the sets $\{u_D>{\psi}\},\{u_N>{\psi}\}$ have positive measure as
they coincide with ${\kappa}$, but $P[u_N-{\psi}]=P[u_D-{\psi}]=\emptyset$.
\end{Remark}

\begin{Remark}
\label{R:true}
Sufficient conditions in order to have that $P[u_N-\psi]$ is not empty can be easily obtained. For instance,
if $\psi$ vanishes on an open set $B\subset\Omega$, then $B\subseteq P[u_N-\psi]$, since $u_N$ is positive by 
the strong maximum principle in \cite{CDDS}. If $\psi$ is continuous on $\overline\Omega$ then $u_N$ is continuous as well by Theorem \ref{T:2}. 
Thus either $u_N\equiv \psi$, or $P[u_N-\psi]=\{u_N>\psi\}$ is not empty.
\end{Remark}

\end{document}